\documentclass[12pt]{article}
\usepackage{amsmath, amsthm, amscd, amsfonts, amssymb, graphicx, color}
\usepackage{epstopdf}
\usepackage{subfigure}
\usepackage{setspace}
\usepackage[left=1in, right=1in, bottom=1.4in]{geometry}
\usepackage{bm}
\usepackage{cite}
\newtheorem{theorem}{Theorem}[section]

\newtheorem{lemma}[theorem]{Lemma}

\newtheorem{proposition}[theorem]{Proposition}
\newtheorem{remark}{Remark}
\newtheorem{definition}{Definition}

\title{ \bf Boundary behavior of $\alpha$-harmonic functions and their Riesz-Fej\'{e}r  inequalities$^{\ast}$ }
\date{}
 \begin{document}

\maketitle
\renewcommand{\thefootnote}{\fnsymbol{footnote}}
\noindent{\footnotesize\rm{ \begin{center}\bf Bo-Yong Long  \end{center}\
 \begin{center}  School of Mathematical Sciences, Anhui University, Hefei  230601, China
\end{center}}
\footnotetext{\hspace*{-5mm}
\begin{tabular}{@{}r@{}p{16.0cm}@{}}
&$^{*}$Supported by NSFC (No.12271001), Natural Science Foundation of Anhui Province (2308085MA03), and Excellent University Research and Innovation Team in Anhui Province (2024AH010002),
China.
\\
&E-mail:  boyonglong@163.com

\end{tabular}}

\begin{quote}
{\small \noindent {\bf Abstract}:\   The solutions of a kind of second-order homogeneous  partial differential  equation  are called (real kernel) $\alpha$-harmonic functions.
 In this paper, the boundary correspondence and boundary behavior of $\alpha$-harmonic functions are studied, and the corresponding Dirichlet problem is solved.
As one of its applications, an asymptotic optimal Riesz-Fej\'{e}r  inequality for $\alpha$-harmonic functions is obtained. In addition, the subharmonic properties of $\alpha$-harmonic functions is explored and an optimal radius is obtained.
}\\
  {{\small \noindent{\bf Keywords}: $\alpha$-harmonic functions;  boundary behavior;  Dirichlet problem; subharmonicity, Riesz-Fej\'{e}r inequality }
 }
\end{quote}

\section{Introduction and main results}
\setcounter{equation}{0}
\hspace{2mm}

{\bf 1.0  $\alpha$-harmonic functions}

Let $\mathbb{D}=\{z:|z|<1\}$ be the  open unit disk  and $\mathbb{T}=\{z:|z|=1\}$ unit circle.
For $\alpha\in\mathbb{R}$ and  $z\in \mathbb{D}$, let

$$T_{\alpha}=-\frac{\alpha^{2}}{4}(1-|z|^{2})^{-\alpha-1}+\frac{\alpha}{2}(1-|z|^{2})^{-\alpha-1}\left(z\frac{\partial}{\partial z}
+\bar{z}\frac{\partial}{\partial\bar{z}}\right)+\frac{1}{4}(1-|z|^{2})^{-\alpha}\triangle$$
be a second-order elliptic  partial differential operator, where $\triangle$ is the usual complex Laplacian operator
$$\triangle:=4\frac{\partial^{2}}{\partial z \partial \bar{z}}=\frac{\partial^{2}}{\partial x^{2}}+\frac{\partial^{2}}{\partial y^{2}}. $$
The corresponding homogeneous differential equation is
\begin{align}\label{1.1}
T_{\alpha}(u)=0 \quad \mbox {in}  \,\,\mathbb{D}
\end{align}
and the associated the Dirichlet boundary value problem is as follows
\begin{equation} \label{1.2}\;  \left\{
\begin{array}{rr}
T_{\alpha}(u)=0 &\quad \mbox {in}  \,\,\mathbb{D}, \\
u=f &\quad \mbox {on}  \,\,\mathbb{T}.
\end{array}\right.
\end{equation}
Here, the boundary data $f\in\mathfrak{D}'(\mathbb{T})$
is a distribution on the boundary $\mathbb{T}$ of $\mathbb{D}$, and the boundary condition in (\ref{1.2}) is interpreted in the distributional sense that $u_{r}\rightarrow f$ in $\mathfrak{D}'(\mathbb{T})$
as $r\rightarrow1^{-}$, where
\begin{equation*}
u_{r}(e^{i\theta})=u(re^{i\theta}), \quad e^{i\theta}\in\mathbb{T},
\end{equation*}
for $r\in[0,1)$. By the results of  \cite{Olofsson2014}, we know that if $\alpha>-1$ and $u\in\mathcal{C}^{2}(\mathbb{D})$ satisfy (\ref{1.2}), then  the function $u$ has the form of Poisson type integral
\begin{equation}\label{1.3}
u(z):=P_{\alpha}[f](z)=\frac{1}{2\pi}\int_{0}^{2\pi}K_{\alpha}(ze^{-it)})f(e^{it})dt, \quad \mbox{for}\,\, z\in\mathbb{D},
\end{equation}
where
 \begin{equation}\label{1.4}
 K_{\alpha}(z)=c_{\alpha}\frac{(1-|z|^{2})^{\alpha+1}}{|1-z|^{\alpha+2}},
 \end{equation}
$c_{\alpha}=\Gamma^{2}(\alpha/2+1)/\Gamma(1+\alpha)$ and $\Gamma(\cdot)$  is the  Gamma function.

Regarding the expression of the solutions of (\ref{1.1}),  the following results have been obtained in \cite{Olofsson2014}.

{\bf Theorem A}\cite{Olofsson2014}
 Let $\alpha\in\mathbb{R}$ and $u\in\mathcal{C}^{2}(\mathbb{D})$. Then $u$ satisfies (\ref{1.1}) if and only if it has a series expansion
 of the form
 \begin{equation}\label{1.5}
 u(z)=\sum_{n=0}^{\infty}c_{n}F(-\frac{\alpha}{2},n-\frac{\alpha}{2}; n+1; |z|^{2})z^{n}+\sum_{n=1}^{\infty}c_{-n}F(-\frac{\alpha}{2},n-\frac{\alpha}{2}; n+1; |z|^{2})\bar{z}^{n}, \quad z\in \mathbb{D},
 \end{equation}
 for some sequence $\{c_{n}\}_{-\infty}^{\infty}$ of complex number satisfying
 \begin{equation}\label{1.6}
 \lim_{|n|\rightarrow\infty}\sup|c_{n}|^{\frac{1}{|n|}}\leq 1.
 \end{equation}
 In particular, the expansion (\ref{1.5}), subject to (\ref{1.6}), converges in $\mathcal{C}^{\infty}(\mathbb{D})$, and every solution $u$ of (\ref{1.1}) is $\mathcal{C}^{\infty}-$smooth in $\mathbb{D}$. Here,  $F(a,b;c;x)$ denote the Gauss hypergeometric functions.

It was shown in \cite{Olofsson2014} that if $\alpha\leq -1$, $u\in\mathcal{C}^{2}(\mathbb{D})$ satisfy (\ref{1.1}) and  the boundary limit  $f=\lim_{r\rightarrow 1^{-}}u_{r}$ exists in $\mathfrak{D}'(\mathbb{T})$, then
  $u(z)=0$ for all $z\in\mathbb{D}$.  That is to say, the case of $\alpha\leq -1$ is trivial. So, in the following we will only focus on the case of $\alpha>-1$.
We call the solutions of (\ref{1.1})  with $\alpha>-1$
 real kernel $\alpha$-harmonic functions \cite{Long2022Proc, Long2021Filomat, Long2024}, or  $\alpha$-harmonic   functions.

On the one hand, if we take $\alpha=2(p-1)$, then $u$ is polyharmonic (or $p-$harmonic), where $p\in\{1,2,...\}$, cf.\cite{Du2012, Chens2021, Lip2022, Lium2024,Iwaniec2019}.
In particular, if $\alpha=0$, then $u$ is harmonic, cf.\cite{Kalaj2015,Kalaj2019,Bshouty2018,Duren2004}.  Thus, $u$ is a kind of generalization of classical harmonic mappings. On the other hand, $\alpha$-harmonic functions themselves have also been generalized to so-called $(\alpha, \beta)$-harmonic functions or $(p, q)$-harmonic functions \cite{Klintborg2021, Arsenovic2024}.

{\bf 1.1  Boundary corresponding of the $\alpha$-harmonic functions  and the Dirichlet problem}

Similar to harmonic functions, when the boundary function on the unit circumference is continuous, the $\alpha$-harmonic function on the unit disk is also continuous at the boundary. Therefore, in this case, there is a good boundary correspondence and boundary extension phenomenon. Actually, we have the following results:

\begin{theorem} \label{Th1.1}
Let  function $f$ be piecewise continuous on the unit circumference $\mathbb{T}$. Then for any $\alpha>-1$, $u(z)=P_{\alpha}[f](z)$  is a bounded $\alpha$-harmonic function on the unit disk $\mathbb{D}$, and at the continuous point $e^{it_{0}}$ of $f$, it holds that
\begin{align}\label{1.7}
\lim_{z\rightarrow e^{it_{0}}}u(z)=u(e^{it_{0}}),\quad z\in\mathbb{D}.
\end{align}
\end{theorem}

Based on the results of Theorem \ref{Th1.1}, regarding the classical Dirichlet problem, the $\alpha$-harmonic functions also have similar conclusions to harmonic functions. Specifically, we have
\begin{theorem}\label{Th1.2}
Let $\alpha>-1$ and $f$ be  continuous on the unit circumference $\mathbb{T}$. Then there exists a unique solution  $u=P_{\alpha}[f]$ to the the system (\ref{1.2})
  such that $u$ is continuous on $\overline{\mathbb{D}}$.
\end{theorem}

\begin{remark}
 Let
\begin{align*}  f(e^{i\theta})&=\sum^{\infty}_{n=0}c_{n}e^{in\theta}+\sum^{\infty}_{n=1} c_{-n}e^{-in\theta},\\
U(z)&=\sum^{\infty}_{n=0}c_{n}z^{n}+\sum^{\infty}_{n=1} c_{-n}\bar{z}^{n},
\end{align*} and
 \begin{align*} u(z)=&\sum^{\infty}_{n=0}\frac{c_{n}}{F(-\frac{\alpha}{2},n-\frac{\alpha}{2}; n+1; 1)}F(-\frac{\alpha}{2},n-\frac{\alpha}{2}; n+1; |z|^{2})z^{n}\\
    &+\sum^{\infty}_{n=1} \frac{c_{-n}}{F(-\frac{\alpha}{2},n-\frac{\alpha}{2}; n+1; 1)}F(-\frac{\alpha}{2},n-\frac{\alpha}{2}; n+1; |z|^{2})\bar{z}^{n}.
    \end{align*}

 (1)\quad Obviously, $U(z)$ and $u(z)$  have the expression forms of harmonic function and $\alpha$-harmonic function,  respectively. $f(e^{i\theta})$ is the boundary function of $U(z)$ and $u(z)$.

 From the proof process of the theorem \ref{Th1.2}, it can be further inferred that, among $U(z)$, $u(z)$, and $f(e^{i\theta})$ these three, there is a one-to-one correspondence relationship between any two of them, regardless of whether the boundary function  $f(e^{i\theta})$ is continuous or not.

(2) \quad If  $f(e^{i\theta})$ is not a constant and continuous on  $\mathbb{T}$, then both image domain of $U(z)$ and $u(z)$ fill the area besieged by the closed curve $f(e^{i\theta})$.
   In other words, $U(z)$ and $u(z)$ have exactly the same image domain. This suggests that we can consider using harmonic mappings on unit disk $\mathbb{D}$ to study $\alpha$-harmonic mappings with the same boundary function.
\end{remark}

{\bf 1.2  Boundary behavior of the $\alpha$-harmonic functions }

   For the case where there are jump discontinuities on the boundary, we have the following conclusion:

\begin{theorem}\label{Th1.3} Suppose  $f$ is piecewise continuous with a finite number of jump discontinuities on the unit circumference $\mathbb{T}$, $e^{i\theta_{0}}$ is a jump discontinuity whose $\lim_{\theta\rightarrow\theta_{0}^{+}}f(e^{i\theta})$ and $\lim_{\theta\rightarrow\theta_{0}^{-}}f(e^{i\theta})$ exist but  $\lim_{\theta\rightarrow\theta_{0}^{+}}f(e^{i\theta})\neq\lim_{\theta\rightarrow\theta_{0}^{-}}f(e^{i\theta})$. If  the point $z$ in the unit disk approaches the boundary point $e^{i\theta_{0}}$ along a linear segment at an angle $\gamma, \, (0<\gamma<\pi)$ with the tangent line, then the $\alpha$-harmonic function $u(z)=P_{\alpha}[f](z)$ tends to the corresponding weighted average
\begin{align*}\frac{\gamma}{\pi}\lim_{\theta\rightarrow\theta_{0}^{-}}f(e^{i\theta})
+(1-\frac{\gamma}{\pi})\lim_{\theta\rightarrow\theta_{0}^{+}}f(e^{i\theta}).
\end{align*}
\end{theorem}

When $\alpha=0$, Theorem \ref{Th1.3} reduces to the case of harmonic mappings, see section 1.4 of \cite{Duren2004}. For more conclusions on the boundary behavior of harmonic mappings, please refer to \cite{Bshouty2012,Bshouty2022}.

{\bf 1.3  Subharmonicity of the $\alpha$-harmonic functions }

The proof of the uniqueness of the Dirichlet problem for harmonic functions often uses the extremum principle of harmonic functions or the properties of subharmonic functions. But in the proof of the uniqueness of Theorem \ref{Th1.2}, we did not use a similar method. Because the real valued $\alpha$-harmonic function may not necessarily be subharmonic on the unit disk.
In fact, we have the following interesting conclusion for real-valued $\alpha$-harmonic functions:

\begin{theorem}\label{Th1.4}
Suppose $f(e^{it})\geq 0$ for all $t\in[0,2\pi]$. Then the $\alpha$-harmonic function $u(z)=P_{\alpha}[f](z)$ is non-negative and  subharmonic on the disk $\mathbb{D}_{r_{0}}=\{z| |z|<r_{0}\}$ with
\begin{equation} \label{1.8}\; r_{0}= \left\{
\begin{array}{rr}
\frac{\sqrt{1+\alpha}-1}{\sqrt{1+\alpha}+1}, &\quad \alpha>0,\quad\\
1,\qquad &\quad \alpha=0,\quad\\
\frac{1-\sqrt{1+\alpha}}{\sqrt{1+\alpha}+1},  &\quad \alpha\in(-1,0).
\end{array}\right.
\end{equation}Furthermore, the constant $r_{0}$ is optimal.
\end{theorem}

Similar to the concepts of star radius and convex radius in geometric function theory, we can call the constant $r_{0}$ in Theorem \ref{Th1.4} the subharmonic radius. The concept of subharmonic radius seems to have not been proposed before.

\begin{proposition}
There does exist that $\alpha$-harmonic functions  are subharmonic over the entire unit disk $\mathbb{D}$.
\end{proposition}

{\bf  1.4 Riesz-Fej\'{e}r inequality for $\alpha$-harmonic Hardy spaces}

In 1921,  Fej\'{e}r and  Riesz have shown  that for any conformal map of the unit circle onto a plane domain, the length of the map of the circumference is at least twice that of the map of any diameter, cf. Theorem 3.13 of \cite{Duren1970}. This result was generalized in several directions, cf \cite{Beckenbach1938,  Wulan2011, Sagher1977, Chensl2023, Melentijevic2023PA, Melentijevic2024MA}. There are also conclusions about harmonic functions in this regard\cite{Melentijevic2021, Kayumov2020,Das2022}. Among them, the conclusions we focus on specifically are

{\bf Theorem B} \cite{Melentijevic2021}\quad For all $1<p<\infty$ and $f\in h^{p}(\mathbb{D})$ and any $s\in[0,2\pi]$, we have:
\begin{align}\label{1.9}
\int^{1}_{-1}|f(re^{is})|^{p}dr\leq  \frac{1}{2}\sec^{p}(\frac{\pi}{2p})\int^{2\pi}_{0}|f(e^{i\theta})|^{p}d\theta.
\end{align}

Here, $h^{p}(\mathbb{D})$ denote the Hardy space composed of complex valued harmonic functions on the unit disk $\mathbb{D}$. Inequality (\ref{1.9}) is optimal, that is, the constant  $\frac{1}{2}\sec^{p}(\frac{\pi}{2p})$ is sharp. In \cite{Kayumov2020}, it is proven that  inequality (\ref{1.9}) is optimal for $p\in(1,2]$. It is worth pointing out that the method of proving  inequality (\ref{1.9}) is optimal for the case $p\in(1,2]$ is fully applicable to the case $p\in(2,\infty)$.

Let $h_{\alpha}^{p}(\mathbb{D})$ denote the Hardy space composed of $\alpha$-harmonic functions on  unit disk $\mathbb{D}$. We have generalized the conclusion of Theorem B to the case of  $\alpha$-harmonic functions. Specifically, we have

\begin{theorem}\label{Th1.5}
Suppose that $u\in h_{\alpha}^{p}(\mathbb{D})$ with $p>1$. If $-1<\alpha\leq 0$ and $\alpha+\frac{2}{p}>0$, then for any $s\in[0,2\pi]$, it holds that
\begin{align}\label{1.10}
\int^{1}_{-1}|u(re^{is})|^{p}dr\leq   C(\alpha, p)\int^{2\pi}_{0}|u(e^{i\theta})|^{p}d\theta,
\end{align} where
\begin{align*}C(\alpha, p)=\frac{2^{\alpha-1}c_{\alpha}}{\pi}B\left(\frac{1}{2}(1+\alpha+\frac{1}{p}),\frac{1}{2}(1-\frac{1}{p})\right)
\sec^{p-1}(\frac{\pi}{2p}),
\end{align*} $B(\cdot, \cdot)$ denotes Beta function. The constant $C(\alpha, p)$ is asymptotically sharp as $\alpha\rightarrow 0$.

\end{theorem}

\section{Prliminaries}
\setcounter{equation}{0}
\hspace{2mm}

In this section, we shall recall some necessary terminology and useful known results.

{\bf 2.1 Schur test}

The following version of Schur test can be found in \cite{Howard1990}.
\begin{lemma}\label{lem2.1}
Let $X$ and $Y$ be measure spaces equipped with  nonnegative, $\sigma$-finite measures and let $T$ be an operator from $L^{p}(Y)$ to $L^{p}(X)$ that can be expressed as
\begin{align*}
Tf(x)=\int_{Y}K(x,y)f(y)dy
\end{align*}for some nonnegative function $K(x,y)$. The adjoint operator $T^{*}$ is given by
\begin{align*}
T^{*}f(y)=\int_{X}K(x,y)f(x)dx.
\end{align*}If finding a measurable $h$ finite almost everywhere, such that
\begin{align*}
T^{*}((Th)^{p-1})\leq A_{p}h^{p-1},
\end{align*}a.e. on $Y$, then for all $f\in L^{p}(Y)$, we have
\begin{align*}
\int_{X}|T(f)|^{p}dx\leq  A_{p}\int_{Y}|f|^{p}dy.
\end{align*}

\end{lemma}

{\bf 2.2 Gauss hypergeometric functions}

The  Gauss hypergeometric function is defined by the  series
$$F(a,b;c; x)=\sum_{n=0}^{\infty}\frac{(a)_{n}(b)_{n}}{(c)_{n}}\frac{x^{n}}{n!}$$ for $ |x|<1$,  and by continuation elsewhere,
where $(a)_{0}=1$ and $(a)_{n}=a(a+1)\cdots(a+n-1)$ for $n=1,2,...$ are the Pochhammer symbols. Obviously, for $n=0, 1,2,...$,
$(a)_{n}=\Gamma(a+n)/\Gamma(a)$.

{\bf  Conclusions}\cite{Andrews1999}

1. If $\Re(c-a-b)>0$, then
\begin{equation}\label{2.1}
\lim_{x\rightarrow 1}F(a,b;c; x)=\frac{\Gamma(c)\Gamma(c-a-b)}{\Gamma(c-a)\Gamma(c-b)}.
\end{equation}

2. If $\Re(c-a-b)<0$ and $|x|<1$, then
\begin{equation}\label{2.2}
F(a,b;c; x)=(1-x)^{c-a-b}F(c-a, c-b; c; x).
\end{equation}

3. It holds that\begin{equation}\label{2.3}
\frac{d F(a,b; c; x)}{dx}=\frac{ab}{c}F(a+1, b+1; c+1; x).
\end{equation}

\begin{lemma}\cite{Olofsson2014}\label{lem2.2}
Let $c>0$, $a\leq c$, $b\leq c$ and $ab\leq 0$ $(ab\geq 0)$. Then the function $F(a,b;c;x)$ is decreasing (increasing) on $x\in [0, 1)$.
\end{lemma}

{\bf 2.3  Subharmonic functions}

\begin{definition}  A function $v(z)$ of one complex or two real variables is said to be subharmonic if in any region $v(z)$ is less than or equal to the harmonic function $u(z)$ which coincides with $v(z)$ on the boundary of the region.
\end{definition}

{\bf Conclusions}

1.   If a continuous real-valued function $v(z)$ is subharmonic in region $\Omega$, then for any region $\Omega'\subset\Omega$, $v(z)$ satisfies the maximum principle in $\Omega'$.

2. If $\Delta v\geq0$ in region $\Omega$, then $v$ is subharmonic in $\Omega$.

{\bf 2.4   Hardy space}

Let $f$ be measurable  complex-valued function defined on unit disk $\mathbb{D}$. The integral means of $f$ are defined as follows:
\begin{align*}
M_{p}(r,f)=\left(\frac{1}{2\pi}\int^{2\pi}_{0}|f(re^{i\theta})|^{p}d\theta\right)^{1/p}, \quad 0<p<\infty;
\end{align*}and
\begin{align*}
M_{\infty}(r,f)=ess\sup\limits_{0\leq\theta\leq2\pi}|f(re^{i\theta})|.
\end{align*}
A function $f$ analytic in $\mathbb{D}$ is said to be of class $H^{p}(\mathbb{D})$ if $M_{p}(r,u)$ is bounded.

It is convenient also to define the analogous classes of harmonic functions or $\alpha$-harmonic functions. A function $u$ harmonic or $\alpha$-harmonic in $\mathbb{D}$ is said to be of class $h^{p}(\mathbb{D})$ or $h_{\alpha}^{p}(\mathbb{D})$ if $M_{p}(r,u)$ is bounded, respectively.

{\bf 2.5   A lemma}

\begin{lemma}\cite{Ponnusamy2001, YangZ2015} \label{lem2.5}
Let $r_{n}$ and $s_{n}$\, $(n=0,1,2,...)$ be real numbers, and let the power series $$R(x)=\sum_{n=0}^{\infty}r_{n}x^{n} \quad \mbox{and}\quad S(x)=\sum_{n=0}^{\infty}s_{n}x^{n}$$
be convergent for $|x|<r$, $(r>0)$ with $s_{n}>0$ for all $n$. If the non-constant sequence $\{r_{n}/s_{n}\}$ is  increasing (decreasing) for all $n$, then the function $x\mapsto R(x)/S(x)$ is strictly increasing (resp. decreasing) on $(0,r)$.
\end{lemma}

\section{Proofs}
\setcounter{equation}{0}
\hspace{2mm}

\begin{proof}[\textbf{Proof of Theorem \ref{Th1.1}}]
It is obvious that $u(z)=P_{\alpha}[f](z)$  is an $\alpha$-harmonic function on the unit disk $\mathbb{D}$.

Let \begin{align*}
M_{\alpha}(r)=\frac{1}{2\pi}\int^{2\pi}_{0}K_{\alpha}(re^{\theta})d\theta,
\end{align*}where $K_{\alpha}$ defined by (\ref{1.4}). Then from the proof process of the Theorem 3.1 of \cite{Olofsson2014}, we  obtain that
\begin{align}\label{3.1}
M_{\alpha}(r)=\frac{\Gamma(\alpha/2+1)^{2}}{\Gamma(\alpha+1)}F(-\alpha/2,-\alpha/2; 1;r^{2}).
\end{align}By Lemma \ref{lem2.2} and equation (\ref{2.1}), we have
\begin{align*}
M_{\alpha}(r)\leq M_{\alpha}(1)= 1.
\end{align*}

Because of the boundedness of $f$, there exists  a constant $M$ such that $|f(e^{it})|\leq M$ for all $t\in [0,2\pi]$.  It follows that
\begin{align}\label{3.2}
|u(z)|&\leq M\frac{c_{\alpha}}{2\pi}\int^{2\pi}_{0}\frac{(1-r^{2})^{\alpha+1}}{(1+r^{2}-2r\cos(\theta-t))^{1+\alpha/2}}dt\nonumber\\
&=M\frac{c_{\alpha}}{2\pi}\int^{2\pi}_{0}\frac{(1-r^{2})^{\alpha+1}}{(1+r^{2}-2r\cos t)^{1+\alpha/2}}dt\nonumber\\
&=M \cdot M_{\alpha}(r)\leq M.
\end{align}That is to say  $u(z)$ is bounded.

Now, we need to prove (\ref{1.7}). Let $z=re^{i\theta}$. Observe that
 \begin{align}\label{3.3}
 &\quad |u(z)-u(e^{it_{0}})|\nonumber\\
&=\left|\int^{2\pi}_{0}K_{\alpha}(re^{i(\theta-t)})f(e^{it})\frac{dt}{2\pi}
-f(e^{it_{0}})\frac{1}{M_{\alpha}(r)}\int^{2\pi}_{0}K_{\alpha}(re^{it})\frac{dt}{2\pi}\right|\nonumber\\
&=\left|\int^{2\pi}_{0}K_{\alpha}(re^{i(\theta-t)})f(e^{it})\frac{dt}{2\pi}
-\frac{1}{M_{\alpha}(r)}\int^{2\pi}_{0}K_{\alpha}(re^{i(\theta-t)})f(e^{it_{0}})\frac{dt}{2\pi}\right|\nonumber\\
&\leq\left|\int^{2\pi}_{0}K_{\alpha}(re^{i(\theta-t)})f(e^{it})\frac{dt}{2\pi}
-\frac{1}{M_{\alpha}(r)}\int^{2\pi}_{0}K_{\alpha}(re^{i(\theta-t)})f(e^{it})\frac{dt}{2\pi}\nonumber\right|\\
&\quad\, \left|\frac{1}{M_{\alpha}(r)}\int^{2\pi}_{0}K_{\alpha}(re^{i(\theta-t)})f(e^{it})\frac{dt}{2\pi}
-\frac{1}{M_{\alpha}(r)}\int^{2\pi}_{0}K_{\alpha}(re^{i(\theta-t)})f(e^{it_{0}})\frac{dt}{2\pi}\right|\nonumber\\
&\leq\left(1-\frac{1}{M_{\alpha}(r)}\right)\left|\int^{2\pi}_{0}K_{\alpha}(re^{i(\theta-t)})f(e^{it})\frac{dt}{2\pi}\right|
+\frac{1}{M_{\alpha}(r)}\left|\int^{2\pi}_{0}K_{\alpha}(re^{i(\theta-t)})(f(e^{it})-f(e^{it_{0}}))\frac{dt}{2\pi}\right|\nonumber\\
&:=I+II
\end{align}

Firstly, noting that
\begin{align*}
\lim_{r\rightarrow 1^{-}}(1-\frac{1}{M_{\alpha}(r)})=0,
\end{align*}and considering (\ref{3.2}), we have
 \begin{align}\label{3.4}
 \lim_{r\rightarrow 1^{-}}I=\lim_{r\rightarrow 1^{-}}\left(1-\frac{1}{M_{\alpha}(r)}\right)|u(z)|=0.
 \end{align}

Secondly, since $f$ is continuous at $e^{it_{0}}$, there exists a $\delta$ such that for a given arbitrarily small positive $\epsilon$ we have $|f(e^{it})-f(e^{it_{0}})|<\epsilon$, provided  $t-t_{0}|<\delta$. It follows that
\begin{align}\label{3.5}
&\quad\frac{1}{M_{\alpha}(r)}\left|\int^{t_{0}+\delta}_{t_{0}-\delta}K_{\alpha}(re^{i(\theta-t)})(f(e^{it})-f(e^{it_{0}}))\frac{dt}{2\pi}\right|\nonumber\\
&\leq\epsilon\frac{1}{M_{\alpha}(r)}\int^{t_{0}+\delta}_{t_{0}-\delta}K_{\alpha}(re^{i(\theta-t)})
\frac{dt}{2\pi}\nonumber\\
&<\epsilon\frac{1}{M_{\alpha}(r)}\int^{2\pi}_{0}K_{\alpha}(re^{i(\theta-t)})\frac{dt}{2\pi}\nonumber\\
&=\epsilon.
\end{align}
Furthermore, since the point $z=re^{i\theta}:=(r, \theta)$ approaches the point $e^{it_{0}}:=(1, t_{0})$, we shall ultimately have $|\theta-t_{0}|<\delta/2$. For such values of $\theta$, and values of $t$ for which $|t-t_{0}|>\delta$, we have
\begin{align*}|\theta-t|\geq |t-t_{0}|-|t_{0}-\theta|\geq \delta-\delta/2=\delta/2,
\end{align*}whence $\cos (\theta-t)\leq\cos(\delta/2)$ and
\begin{align*}1+r^{2}-2r\cos(\theta-t)\geq 1+r^{2}-2r\cos(\delta/2)= 1+r^{2}-2r(1-2\sin^{2}(\delta/4))>4r\sin^{2}(\delta/4).
\end{align*}
Let $A=4r\sin^{2}(\delta/4)$.
We therefore obtain
\begin{align}\label{3.6}
&\quad\frac{1}{M_{\alpha}(r)}\left|\int^{t_{0}-\delta+2\pi}_{t_{0}+\delta}K_{\alpha}(re^{(\theta-t)})(f(e^{it})-f(e^{it_{0}}))\frac{dt}{2\pi}\right|\nonumber\\
&\leq\frac{2M}{M_{\alpha}(r)}\int^{t_{0}-\delta+2\pi}_{t_{0}+\delta}\frac{(1-r^{2})^{\alpha+1}}{(1+r^{2}-2\cos(\theta-t))^{\alpha/2+1}}\frac{dt}{2\pi}\nonumber\\
&\leq \frac{2M}{M_{\alpha}(r)}\frac{(1-r^{2})^{\alpha+1}}{A}\longrightarrow 0,  \quad  r\rightarrow 1.
\end{align}
Combining  (\ref{3.5}) and (\ref{3.6}), it yields that if $z\rightarrow e^{it_{0}}$, that is, if $r\rightarrow 1$ and $\theta\rightarrow t_{0}$, then  $II\rightarrow0$.

 Therefore, we have if $z\rightarrow e^{it_{0}}$, then $|u(z)-u(e^{it_{0}})|\rightarrow 0$. The theorem is proved.
\end{proof}

\begin{proof}[\textbf{Proof of Theorem \ref{Th1.2}}]

Firstly,  Let $u(z)=P_{\alpha}[f](z)$, that is to say, $u(z)$ is the $\alpha$-harmonic extension of $f$. Then it follows from \cite{Olofsson2014} that  $T_{\alpha}(u)=0 $ in $\mathbb{D}$ and $u\in C^{\infty}(\mathbb{D})$.

   Secondly,  Theorem \ref{Th1.1} implies that  $u(z)=P_{\alpha}[f](z)$ is continuous on $\overline{\mathbb{D}}$.

  At last, prove the uniqueness of solution. On the one hand, for given $f$, we consider the harmonic extension $U(z):=P_{0}[f](z)$. It is well known that the harmonic extension is unique. Furthermore, the harmonic function $U(z)$ has expression such as \begin{align*} U(z)=\sum^{\infty}_{n=0}a_{n}z^{n}+\sum^{\infty}_{n=1}a_{-n}\bar{z}^{n} \end{align*}with
   \begin{align*} f(e^{i\theta})=\sum^{\infty}_{n=0}a_{n}e^{i n\theta}+\sum^{\infty}_{n=0}b_{n}e^{-i n\theta} \end{align*}as its boundary function.  On the other hand, by Theorem A, any $\alpha$-harmonic function $\widetilde{u}(z)$ on $\mathbb{D}$ must has expression such as  $$\widetilde{u}(z)=\sum^{\infty}_{n=0}c_{n}F(-\frac{\alpha}{2},n-\frac{\alpha}{2}; n+1; |z|^{2})z^{n}+\sum^{\infty}_{n=1}c_{-n}F(-\frac{\alpha}{2},n-\frac{\alpha}{2}; n+1; |z|^{2})\bar{z}^{n}$$
   Let $$c_{n}=\frac{a_{n}}{F(-\frac{\alpha}{2},n-\frac{\alpha}{2}; n+1; 1)} \quad\mbox{and } \quad c_{-n}=\frac{a_{-n}}{F(-\frac{\alpha}{2},n-\frac{\alpha}{2}; n+1; 1)},\quad n=1,2,3....$$ Then
    \begin{align*} u(z)=&\sum^{\infty}_{n=0}\frac{a_{n}}{F(-\frac{\alpha}{2},n-\frac{\alpha}{2}; n+1; 1)}F(-\frac{\alpha}{2},n-\frac{\alpha}{2}; n+1; |z|^{2})z^{n}\\
    &+\sum^{\infty}_{n=1} \frac{a_{-n}}{F(-\frac{\alpha}{2},n-\frac{\alpha}{2}; n+1; 1)}F(-\frac{\alpha}{2},n-\frac{\alpha}{2}; n+1; |z|^{2})\bar{z}^{n}
    \end{align*} is an $\alpha$-harmonic function with same boundary date $f$. So, $u(z)=P_{\alpha}[f](z)$ is the $\alpha$-harmonic extension of $f$.  Because  $u(z)$ and $U(z)$ have same boundary date $f$ and $u(z)\mapsto U(z)$ is one-to-one, the uniqueness of $\alpha$-harmonic extension of $f$  comes from the uniqueness of harmonic extension of $f$.
\end{proof}

\begin{proof}[\textbf{Proof of Theorem \ref{Th1.3}}]

Suppose $f(e^{i\theta})\in\mathbb{R} $ and $e^{i\theta_{0}}$ is the unique jump discontinuity of $f$ with $$\lim_{\theta\rightarrow\theta_{0}^{+}}f(e^{i\theta})=M \quad\mbox{and } \lim_{\theta\rightarrow\theta_{0}^{-}}f(e^{i\theta})=m.$$

Let
\begin{align*}
G(z)=\arctan\left(\frac{y-\cos\theta_{0}}{x-\sin\theta_{0}}\right), \quad z=x+yi\in\overline{\mathbb{D}}\backslash\{e^{i\theta_{0}}\},
\end{align*}and
 \begin{align*} g(e^{i\theta})=G(z)|_{\mathbb{T}}.
\end{align*}
Through direct verification, we know that function $G(z)$ is  harmonic on $\mathbb{D}$.
Furthermore,
since $(y-\cos\theta_{0})(x-\sin\theta_{0})^{-1}$ is the slope of the straight line connecting the points $z=x+yi$ and $e^{\theta_{0}}$, the value of $G(z)$ coincides with that of the angle which this line forms with the positive axis.  Observe that if the point $z=e^{i\theta}:=(\cos\theta, \sin\theta):=(x,y)$ passes through point $e^{i\theta_{0}}$ while describing the unite circle $\mathbb{T}$, then this angle clearly jumps by the amount $\pi$.  So, if let
 \begin{align*} \lim_{\theta\rightarrow\theta_{0}^{+}}g(e^{i\theta})=\beta_{0},\end{align*}
 then
 \begin{align*} \lim_{\theta\rightarrow\theta_{0}^{-}}g(e^{i\theta})=\beta_{0}+\pi.\end{align*}
Construct auxiliary function
\begin{align*}
h(e^{i\theta})=f(e^{i\theta})+\frac{M-m}{\pi}g(e^{i\theta}).
\end{align*} It follows that
\begin{align*}
 &\lim_{\theta\rightarrow\theta_{0}^{+}}h(e^{i\theta})=M+\frac{M-m}{\pi}\beta_{0}=(1+\frac{\beta_{0}}{\pi})M-\frac{\beta_{0}}{\pi}m.\\
 &\lim_{\theta\rightarrow\theta_{0}^{-}}h(e^{i\theta})=m+\frac{M-m}{\pi}(\beta_{0}+\pi)=(1+\frac{\beta_{0}}{\pi})M-\frac{\beta_{0}}{\pi}m.
\end{align*}Thus,
\begin{align*}\lim_{\theta\rightarrow\theta_{0}^{+}}h(e^{i\theta})
=\lim_{\theta\rightarrow\theta_{0}^{-}}h(e^{i\theta}):=h(e^{i\theta_{0}}).
 \end{align*} Therefore, $h(e^{i\theta})$ is continuous on $e^{i\theta_{0}}$. So does on $\mathbb{T}$.
Using function $h$ as the boundary function, conduct $\alpha$-harmonic continuation, then we have
 \begin{align*}
P_{\alpha}[h](z)=P_{\alpha}[f](z)+\frac{M-m}{\pi}P_{\alpha}[g](z).
\end{align*} Now we denote by $L_{\gamma}$ the linear segment which forms an angle $\gamma, \,(0<\gamma<\pi)$ with the tangent line  of unit circumference $\mathbb{T}$ at
  point $e^{i\theta_{0}}$.
It follows that \begin{align*}
\lim\limits_{\substack{z\rightarrow e^{i\theta_{0}}\\z\in L_{\gamma}\bigcap\mathbb{D}}}P_{\alpha}[f](z)=\lim\limits_{\substack{z\rightarrow e^{i\theta_{0}}\\z\in L_{\gamma}\bigcap\mathbb{D}}}P_{\alpha}[h](z)-\frac{M-m}{\pi}\lim\limits_{\substack{z\rightarrow e^{i\theta_{0}}\\z\in L_{\gamma}\bigcap\mathbb{D}}}P_{\alpha}[g](z).
\end{align*}
On the one hand, because $h(e^{i\theta})$ is continuous on $\mathbb{T}$, considering Theorem \ref{1.1}, we have that\begin{align*} \lim\limits_{\substack{z\rightarrow e^{i\theta_{0}}\\z\in L_{\gamma}\bigcap\mathbb{D}}}P_{\alpha}[h](z)=h(e^{i\theta_{0}})=(1+\frac{\beta_{0}}{\pi})M-\frac{\beta_{0}}{\pi}m.
\end{align*}
On the other hand, let $\bigcup\limits^{\circ}(e^{i\theta_{0}},\delta)=\{z| \,0<|z-e^{i\theta_{0}}|<\delta\}$, then for any $\alpha>-1$ and $e^{i\theta}\in \bigcup\limits^{\circ}(e^{i\theta_{0}}, \delta)\bigcap\mathbb{T}$,
$e^{i\theta}$ is the continuous point of $g$. Then  by Theorem \ref{1.1} again, we have that
\begin{align*} \lim\limits_{\substack{z\rightarrow e^{i\theta}\\z\in \mathbb{D}}}P_{\alpha}[g](z)=g(e^{i\theta}).
\end{align*}Specifically, when $\alpha=0$, we have
\begin{align*} \lim\limits_{\substack{z\rightarrow e^{i\theta}\\z\in \mathbb{D}}}P_{0}[g](z)=g(e^{i\theta}).
\end{align*}Note that $G(z)$ is harmonic, and its boundary function is just $g(e^{i\theta})$. So
\begin{align*}
P_{0}[g](z)=G(z), \quad z\in\overline{\mathbb{D}}\backslash\{e^{i\theta_{0}}\}.
\end{align*}Therefore, $\forall \, \epsilon>0, \exists \,\delta>0$, such that   $\forall  \,z\in \bigcup\limits^{\circ}(e^{i\theta_{0}},\delta)\bigcap\overline{\mathbb{D}}$, it hold that \begin{align*}|P_{\alpha}[g](z)-P_{0}[g](z)|=|P_{\alpha}[g](z)-G(z)|\leq|P_{\alpha}[g](z)-g(e^{i\theta})|
+|g(e^{i\theta})-G(z)|<2\epsilon.
\end{align*}where  $e^{i\theta}\in \bigcup\limits^{\circ}(e^{i\theta_{0}}, \delta)\bigcap\mathbb{T}$.
Because $\epsilon$ is arbitrarily small, we have
\begin{align*} \lim\limits_{\substack{z\rightarrow e^{i\theta_{0}}\\z\in \mathbb{D}}}P_{\alpha}[g](z)=\lim\limits_{\substack{z\rightarrow e^{i\theta_{0}}\\z\in \mathbb{D}}}G(z).
\end{align*}Specifically, we have
\begin{align*} \lim\limits_{\substack{z\rightarrow e^{i\theta_{0}}\\z\in L_{\gamma}\bigcap\mathbb{D}}}P_{\alpha}[g](z)=\lim\limits_{\substack{z\rightarrow e^{i\theta_{0}}\\z\in L_{\gamma}\bigcap\mathbb{D}}}G(z)=\beta_{0}+\gamma.
\end{align*}
Theorem, \begin{align*}\lim\limits_{\substack{z\rightarrow e^{i\theta_{0}}\\z\in L_{\gamma}\bigcap\mathbb{D}}}P_{\alpha}[f](z)
=(1+\frac{\beta_{0}}{\pi})M-\frac{\beta_{0}}{\pi}m-\frac{M-m}{\pi}(\beta_{0}+\gamma)=M(1-\frac{\gamma}{\pi})+m\frac{\gamma}{\pi}.
\end{align*}

If $f(e^{i\theta})\in \mathbb{C}$, then divide complex number $f(e^{i\theta})$ into real and imaginary parts, and repeat the above proof process for both parts,  separately.

If $f(e^{i\theta}) $  have finite number jump discontinuities
$e^{i\theta_{k}}, \, k=1,2,...N$
with $\lim_{\theta\rightarrow\theta_{k}^{+}}f(e^{i\theta})=M_{k}$ and $\lim_{\theta\rightarrow\theta_{k}^{-}}f(e^{i\theta})=m_{k}$.
Then let
\begin{align*}
G_{k}(z)=\arctan\left(\frac{y-\cos\theta_{k}}{x-\sin\theta_{k}}\right), \quad z=x+yi\in\overline{\mathbb{D}}\backslash\{e^{i\theta_{k}}\}, \quad  k=1,2,...N.
\end{align*}
and construct auxiliary function
\begin{align*}
h(e^{i\theta})=f(e^{i\theta})+\sum_{k=1}^{k=N}\frac{M_{k}-m_{k}}{\pi}g_{k}(e^{i\theta}),
\end{align*}where $g_{k}=G_{k}|_{\mathbb{T}}$.
For every fixed jump discontinuity $e^{i\theta_{k}}$, repeat the proof steps above. Then  the proof  is complete.
\end{proof}

\begin{proof}[\textbf{Proof of Theorem \ref{Th1.4}}]

Let $z=re^{i\theta}$ and $K_{\alpha}$ defined as in (\ref{1.4}). Noting that $K_{\alpha}$  is non-negative, if $f\geq 0$, then $u(z)=P_{\alpha}[f](z)$ is  non-negative and $\alpha$-harmonic.

In  Theorem 1.1 of \cite{Olofsson2014}, it is  actually proved that
\begin{align*} T_{\alpha}K_{\alpha}(z)=0.
\end{align*}It follows that
 \begin{align} \label{3.7}
 4(1-r^{2})\Delta K_{\alpha}(z)=(\alpha^{2}I-2\alpha(z\partial+\bar{z}\bar{\partial}))K_{\alpha}(z),
\end{align}where $I$ is identity operator.
Direct computation yields
\begin{align*}\partial K_{\alpha}(z)=\left(-\frac{(\alpha+1)\bar{z}}{1-r^{2}}+\frac{\frac{\alpha}{2}+1}{1-z}\right)K_{\alpha}(z),\\
\bar{\partial } K_{\alpha}(z)=\left(-\frac{(\alpha+1)z}{1-r^{2}}+\frac{\frac{\alpha}{2}+1}{1-\bar{z}}\right)K_{\alpha}(z).
\end{align*}Thus,
\begin{align*}(\alpha^{2}I-2\alpha(z\partial+\bar{z}\bar{\partial}))K_{\alpha}(z)
=\alpha\left(\alpha+4(\alpha+1)\frac{r^{2}}{1-r^{2}}-2(\alpha+2)\Re\frac{z}{1-z}\right)K_{\alpha}(z).
\end{align*}
If $ \alpha>0$, then
\begin{align}\label{3.8} &\quad\alpha\left(\alpha+4(\alpha+1)\frac{r^{2}}{1-r^{2}}-2(\alpha+2)\Re\frac{z}{1-z}\right)K_{\alpha}(z)\nonumber\\
&\geq\alpha\left(\alpha+4(\alpha+1)\frac{r^{2}}{1-r^{2}}-2(\alpha+2) \frac{r}{1-r}\right)K_{\alpha}(z)\nonumber\\
&=\frac{\alpha}{1-r^{2}}(\alpha r^{2}-2(\alpha+2)r+\alpha)K_{\alpha}(z)\nonumber\\
&\geq 0
\end{align}for $r\in[0, r_{0}]$ with $r_{0}=\frac{\sqrt{1+\alpha}-1}{\sqrt{1+\alpha}+1}$.

If $ \alpha=0$, then
\begin{align*} \alpha\left(\alpha+4(\alpha+1)\frac{r^{2}}{1-r^{2}}-2(\alpha+2)\Re\frac{z}{1-z}\right)K_{\alpha}(z)=0
\end{align*}for $r\in[0, r_{0})$ with $r_{0}=1$.

If $ \alpha\in(-1,0)$, then
\begin{align} \label{3.9} &\quad\alpha\left(\alpha+4(\alpha+1)\frac{r^{2}}{1-r^{2}}-2(\alpha+2)\Re\frac{z}{1-z}\right)K_{\alpha}(z)\nonumber\\
&\geq\alpha\left(\alpha+4(\alpha+1)\frac{r^{2}}{1-r^{2}}-2(\alpha+2)\Re\frac{-r}{1+r}\right)K_{\alpha}(z)\nonumber\\
&=\frac{\alpha}{1-r^{2}}(\alpha r^{2}+2(\alpha+2)r+\alpha)K_{\alpha}(z)\nonumber\\
&\geq 0
\end{align}for $r\in[0, r_{0}]$ with $r_{0}=\frac{1-\sqrt{1+\alpha}}{\sqrt{1+\alpha}+1}$.

Considering (\ref{3.7}), whether $\alpha>0$, $\alpha=0$ or $-1<\alpha<0$, we have that $\Delta K_{\alpha}(z)\geq 0$ for $|z|\leq r_{0}$.

So, for fixed $t$, we still have  $\Delta K_{\alpha}(ze^{-it})\geq 0$ for $|z|\leq r_{0}$.  It follows that \begin{align}\label{3.10}
\Delta K_{\alpha}(ze^{-it})\cdot f(e^{it})\geq 0
\end{align}
for $|z|\leq r_{0}$ because  $f(e^{it})\geq 0$. Noting that the operator $\alpha^{2}I-2\alpha(z\partial+\bar{z}\bar{\partial})$ at the right side of (\ref{3.7}) is linear, by integrating both sides of the inequality (\ref{3.10}) simultaneously, we obtain  that $\Delta u(z)=\Delta P_{\alpha}[f](z)\geq 0$ for $|z|\leq r_{0}$.

Let \begin{align*}
f(e^{it})=|\frac{1-ze^{-it}}{1-z}|^{\alpha+2}.
\end{align*} Then $u(z)=P_{\alpha}[f](z)=K_{\alpha}(z)$.  If $\alpha>0$,  taking $z=r$, then the first ``$\geq$'' in (\ref{3.8}) should be ``$=$''. If $ \alpha\in(-1,0)$, taking $z=-r$, then the first ``$\geq$'' in (\ref{3.9}) should be ``$=$''. This fact means   the radius $r_{0}$  is optimal.
\end{proof}

\begin{proof}[\textbf{Proof of Proposition 1.5}]

Taking $n=0$, then $F(-\frac{\alpha}{2},n-\frac{\alpha}{2}; n+1; |z|^{2})z^{n}=F(-\frac{\alpha}{2},-\frac{\alpha}{2}; 1; |z|^{2})$. By Theorem A, $F:=F(-\frac{\alpha}{2},-\frac{\alpha}{2}; 1; |z|^{2})$ is an $\alpha$-harmonic function. Then similar to (\ref{3.7}), we have
\begin{align}\label{3.11}
 4(1-r^{2})\Delta F=(\alpha^{2}I-2\alpha(z\partial+\bar{z}\bar{\partial}))F,
\end{align}where $I$ is identity operator.

Let $t=|z|^{2}=r^{2}$. Then direct computation yields $F_{z}=\bar{z}F_{t}$, $F_{\bar{z}}=zF_{t}$, and
\begin{align}\label{3.12}
(\alpha^{2}I-2\alpha(z\partial+\bar{z}\bar{\partial}))F=\alpha^{2}F-4\alpha tF_{t}=\alpha^{2}F\left(1-\frac{4tF_{t}}{\alpha F}\right)
\end{align}

If $F=F(-\frac{\alpha}{2},-\frac{\alpha}{2}; 1; t):=\sum^{\infty}_{n=0}s_{n}t^{n}$, then $F_{t}=\sum^{\infty}_{n=1}ns_{n}t^{n-1}$ and $tF_{t}=\sum^{\infty}_{n=1}ns_{n}t^{n}:=\sum^{\infty}_{n=1}r_{n}t^{n}$
Then by Lemma \ref{lem2.5},  $r_{n}/s_{n}=n$ is increasing. So, $tF_{t}/F$ is increasing. Considering (\ref{2.1}) and (\ref{2.3}), if $\alpha>0$, then we have
\begin{align} \label{3.13}
\frac{tF_{t}}{F}\leq \left.\frac{tF_{t}}{F}\right|_{t=1}=\frac{(-\frac{\alpha}{2})^{2}F(-\frac{\alpha}{2}+1,-\frac{\alpha}{2}+1;2;1)}
{F(-\frac{\alpha}{2},-\frac{\alpha}{2};1;1)}=\frac{\alpha}{4}.
\end{align}Observe that $F(-\frac{\alpha}{2},-\frac{\alpha}{2}; 1; |z|^{2})>0$. Thus it follows from equations (\ref{3.11})-(\ref{3.13}) that  $F(-\alpha/2,-\alpha/2; 1; |z|^{2})$ is subharmonic on $\mathbb{D}$ for $\alpha>0$.
\end{proof}

\begin{proof}[\textbf{Proof of Theorem \ref{Th1.5}}]

Because of the rotational invariance of norm of functions in $h_{\alpha}^{p}(\mathbb{D})$, or rotational invariance of $\alpha$-harmonic functions, we can consider only the case of $s=0$, without any loss of generality.

We want to apply the Schur test in the following setting.
Let $X=[-1,1]$ with Lebesgue measure and $Y=\partial\mathbb{D}=\mathbb{T}$ with normalised arclength measure. Suppose that an $\alpha$-harmonic function $u\in h_{\alpha}^{p}(\mathbb{D})$, we first get the appropriate $u^{*}(e^{i\theta})\in L^{p}(\mathbb{T})$, defined by its radial limits. By acting with the operator $T$ of Poisson type $\alpha$-harmonic extension, we have
\begin{align*}
Tu^{*}(r)=\int^{2\pi}_{0}c_{\alpha}\frac{(1-r^{2})^{\alpha+1}}
{(1-2r\cos t+r^{2})^{\frac{\alpha+2}{2}}}u^{*}(e^{it})\frac{dt}{2\pi},
\end{align*}where $c_{\alpha}=\Gamma^{2}(\alpha/2+1)/\Gamma(\alpha+1)$. Actually,  by (\ref{1.3}), we have $Tu^{*}(r)=u(r)$. Noting  that $T$ has positive kernel $K_{\alpha}$. It follows that
\begin{align*}
T^{*}u(e^{it})=\int^{1}_{-1}c_{\alpha}\frac{(1-r^{2})^{\alpha+1}}
{(1-2r\cos t+r^{2})^{\frac{\alpha+2}{2}}}u(r)dr.
\end{align*}

Let \begin{align*}h(z)=\Re (1-z^{2})^{-\frac{1}{p}}.
\end{align*} Direct computation yields that if $z=e^{it}\in \mathbb{T}$, then
\begin{align*}
\Re (1-e^{2it})^{-\frac{1}{p}}=\Re (e^{it}(e^{-it}-e^{it}))^{-\frac{1}{p}}=\Re (2\sin t \, e^{i(t-\frac{\pi}{2})})^{-\frac{1}{p}}=2^{-\frac{1}{p}}(\sin t)^{-\frac{1}{p}}\cos\left(\frac{\pi}{2p}-\frac{t}{p}\right),
\end{align*}for $t\in[0,\pi]$;
\begin{align*}
\Re (1-e^{2it})^{-\frac{1}{p}}=\Re (2\sin t e^{i(t-\frac{\pi}{2})})^{-\frac{1}{p}}=2^{-\frac{1}{p}}|\sin t|^{-\frac{1}{p}}\cos\left(\frac{\pi}{2p}-\frac{t-\pi}{p}\right),
\end{align*}for $t\in[\pi,2\pi]$. If $z$ is on the real line, then
\begin{align*}
\Re (1-z^{2})^{-\frac{1}{p}}=(1-r^{2})^{-\frac{1}{p}}.
\end{align*}

Let \begin{align*}V(z)&=\sum_{n=0}^{+\infty}\frac{\Gamma(\frac{\alpha}{2}+1)
\Gamma(\frac{\alpha}{2}+2n+1)}{\Gamma(\alpha+1)(2n)!}(-1)^{n}\tbinom{-\frac{1}{p}}{n}
F(-\frac{\alpha}{2},2n-\frac{\alpha}{2};2n+1; r^{2})z^{2n}\\
&:=\sum_{n=0}^{+\infty}c_{2n}
F(-\frac{\alpha}{2},2n-\frac{\alpha}{2};2n+1; r^{2})z^{2n}.
\end{align*}
Direct computation yields
\begin{align*}\lim_{n\rightarrow +\infty}\left|\frac{c_{2(n+1)}}{c_{2n}}\right|&=\lim_{n\rightarrow +\infty}
\frac{\Gamma(\frac{\alpha}{2}+2(n+1)+1)}{(2(n+1))!}\frac{(\frac{1}{p})_{n+1}}{(n+1)!}
\frac{(2n)!}{\Gamma(\frac{\alpha}{2}+2n+1)}\frac{n!}{(\frac{1}{p})_{n}}\\
&=\lim_{n\rightarrow+\infty}\frac{(\frac{\alpha}{2}+2n+2)(\frac{\alpha}{2}+2n+1)
(\frac{1}{p}+n)}{(2n+2)(2n+1)(n+1)}=1.
\end{align*}Then we have
\begin{align*}\lim_{n\rightarrow +\infty}\sup|c_{2n}|^{\frac{1}{2n}}\leq\lim_{n\rightarrow +\infty}\sup\left|\frac{c_{2(n+1)}}{c_{2n}}\right|=1.
\end{align*}It follows from Theorem A that $V(z)$ is an $\alpha$-harmonic function.

On the one hand, if $-1<\alpha\leq0$, then by Lemma \ref{lem2.2}, $F(-\frac{\alpha}{2},2n-\frac{\alpha}{2};2n+1; r^{2})$ is increasing.  So ,
\begin{align*}
F(-\frac{\alpha}{2},2n-\frac{\alpha}{2};2n+1; r^{2})\leq F(-\frac{\alpha}{2},2n-\frac{\alpha}{2};2n+1;1)
=\frac{\Gamma(\alpha+1)\Gamma(2n+1)}{\Gamma(\frac{\alpha}{2}+1)\Gamma(\frac{\alpha}{2}+2n+1)}.
\end{align*}Furthermore, Observe that $(-1)^{n}\tbinom{-\frac{1}{p}}{n}$ is positive.
Then for $z=r\in(-1, 1)$, we have
\begin{align*}V(r)&=\sum_{n=0}^{+\infty}\frac{\Gamma(\frac{\alpha}{2}+1)
\Gamma(\frac{\alpha}{2}+2n+1)}{\Gamma(\alpha+1)(2n)!}(-1)^{n}\tbinom{-\frac{1}{p}}{n}
F(-\frac{\alpha}{2},2n-\frac{\alpha}{2};2n+1; r^{2})r^{2n}\\
&\leq \sum_{n=0}^{+\infty}\frac{\Gamma(\frac{\alpha}{2}+1)
\Gamma(\frac{\alpha}{2}+2n+1)}{\Gamma(\alpha+1)(2n)!}(-1)^{n}\tbinom{-\frac{1}{p}}{n}
\frac{\Gamma(\alpha+1)\Gamma(2n+1)}{\Gamma(\frac{\alpha}{2}+1)\Gamma(\frac{\alpha}{2}+2n+1)}r^{2n}\\
&=\sum_{n=0}^{+\infty}(-1)^{n}\tbinom{-\frac{1}{p}}{n}r^{2n}=(1-r^{2})^{-\frac{1}{p}}.
\end{align*}
On the other hand,
 \begin{align*}\lim_{r\rightarrow1^{-}}V(re^{i\theta})&=\lim_{r\rightarrow1^{-}}\sum_{n=0}^{+\infty}\frac{\Gamma(\frac{\alpha}{2}+1)
\Gamma(\frac{\alpha}{2}+2n+1)}{\Gamma(\alpha+1)(2n)!}(-1)^{n}\tbinom{-\frac{1}{p}}{n}
F(-\frac{\alpha}{2},2n-\frac{\alpha}{2};2n+1; r^{2})(re^{i\theta})^{2n}\\
&=\sum_{n=0}^{+\infty}\frac{\Gamma(\frac{\alpha}{2}+1)
\Gamma(\frac{\alpha}{2}+2n+1)}{\Gamma(\alpha+1)(2n)!}(-1)^{n}\tbinom{-\frac{1}{p}}{n}
F(-\frac{\alpha}{2},2n-\frac{\alpha}{2};2n+1;1)(e^{i\theta})^{2n}\\
&=(1-e^{2\theta i})^{-\frac{1}{p}}.
\end{align*}That is to say, $(1-e^{2\theta i})^{-\frac{1}{p}}$ is the radial limits of $V(re^{i\theta})$. In other words, $f(e^{i\theta}):=(1-e^{2\theta i})^{-\frac{1}{p}}$ is the boundary function of $V(re^{i\theta})$.

Now, for $-1<r<1$, consider the Remark 1 (1),  we have
\begin{align*}
Th(r)&=\int^{2\pi}_{0}c_{\alpha}\frac{(1-r^{2})^{\alpha+1}}
{(1-2r\cos t+r^{2})^{\frac{\alpha+2}{2}}}\Re (1-e^{2it})^{-\frac{1}{p}}\frac{dt}{2\pi}\\
&=\Re \left(\int^{2\pi}_{0}c_{\alpha}\frac{(1-r^{2})^{\alpha+1}}
{(1-2r\cos t+r^{2})^{\frac{\alpha+2}{2}}}(1-e^{2it})^{-\frac{1}{p}}\frac{dt}{2\pi}\right)\\
&=\Re[ V(re^{i\theta})|_{\theta=0 \,or\,\pi}]\\
&=\Re V(r)\\
&\leq (1-r^{2})^{-\frac{1}{p}}.
\end{align*}Thus,
\begin{align}\label{3.14}
T^{*}((Th)^{p-1})&=\int^{1}_{-1}c_{\alpha}\frac{(1-r^{2})^{\alpha+1}}
{(1-2r\cos t+r^{2})^{\frac{\alpha+2}{2}}}(Th(r))^{p-1}dr\nonumber\\
&\leq \int^{1}_{-1}c_{\alpha}\frac{(1-r^{2})^{\alpha+1}}
{(1-2r\cos t+r^{2})^{\frac{\alpha+2}{2}}}(1-r^{2})^{-\frac{p-1}{p}}dr\nonumber\\
&=\int^{1}_{-1}c_{\alpha}\frac{(1-r^{2})^{\alpha+\frac{1}{p}}}
{(1-2r\cos t+r^{2})^{\frac{\alpha+2}{2}}}dr.
\end{align}

If $t\in[\pi,2\pi]$,  let $t=\theta+\pi$, then $\theta\in[0,\pi]$. Furthermore, let  $r=-s$. Then we have
\begin{align*}
&\quad\int^{1}_{-1}c_{\alpha}\frac{(1-r^{2})^{\alpha+\frac{1}{p}}}{(1-2r\cos t+r^{2})^{\frac{\alpha+2}{2}}}dr=\int^{1}_{-1}c_{\alpha}\frac{(1-r^{2})^{\alpha+\frac{1}{p}}}
{(1+2r\cos \theta+r^{2})^{\frac{\alpha+2}{2}}}dr\\
&=-\int^{-1}_{1}c_{\alpha}\frac{(1-s^{2})^{\alpha+\frac{1}{p}}}
{(1-2s\cos \theta+s^{2})^{\frac{\alpha+2}{2}}}ds=\int^{1}_{-1}c_{\alpha}\frac{(1-s^{2})^{\alpha+\frac{1}{p}}}
{(1-2s\cos \theta+s^{2})^{\frac{\alpha+2}{2}}}ds.\end{align*} So, in the following, we just need to consider the case of  $t\in[0,\pi]$.

Let \begin{align*}\frac{1+r}{1-r}=y\cot\frac{t}{2},\quad t\in[0,\pi].
\end{align*} Then
\begin{align*}r=\frac{y\cot\frac{t}{2}-1}{1+y\cot\frac{t}{2}}, \quad  \text{and}  \quad dr=\frac{2\cot\frac{t}{2}}{(1+y\cot\frac{t}{2})^{2}}dy.
\end{align*} Therefore,
\begin{align*}
&\quad\int^{1}_{-1}c_{\alpha}\frac{(1-r^{2})^{\alpha+\frac{1}{p}}}
{(1-2r\cos t+r^{2})^{\frac{\alpha+2}{2}}}dr\\
&=\int^{+\infty}_{0}c_{\alpha}\frac{\left(1-\left(\frac{y\cot\frac{t}{2}-1}{1+y\cot \frac{t}{2}}\right)^{2}\right)^{\alpha+\frac{1}{p}}}
{\left(1-2\frac{y\cot\frac{t}{2}-1}{1+y\cot\frac{t}{2}}\cos t+\left(\frac{y\cot\frac{t}{2}-1}{1+y\cot \frac{t}{2}}\right)^{2}\right)^{\frac{\alpha+2}{2}}}\frac{2\cot\frac{t}{2}}{(1+y\cot \frac{t}{2})^{2}}dy\\
&=2^{\frac{3\alpha}{2}+\frac{2}{p}}c_{\alpha}\left(\cot\frac{t}{2}\right)^{\alpha+1+\frac{1}{p}}\int^{+\infty}_{0}
\frac{y^{\alpha+\frac{1}{p}}}{\left(y^{2}(1-\cos t)\cot^{2}\frac{t}{2}+1+\cos t\right)^{\frac{\alpha+2}{2}}}
\frac{dy}{(1+y\cot\frac{t}{2})^{\alpha+\frac{2}{p}}}
\end{align*}
\begin{align*}
&=\frac{2^{\alpha+\frac{2}{p}-1}c_{\alpha}\left(\cot\frac{t}{2}\right)^{\alpha+1+\frac{1}{p}}}{(\cos\frac{t}{2})^{\alpha+2}}
\int^{+\infty}_{0}\frac{y^{\alpha+\frac{1}{p}}}{(1+y^{2})^{\frac{\alpha+2}{2}}}
\frac{dy}{(1+y\cot\frac{t}{2})^{\alpha+\frac{2}{p}}}\\
&=2^{\alpha+\frac{2}{p}-1}c_{\alpha}\left(\sin\frac{t}{2}\right)^{-\alpha-1-\frac{1}{p}}
\left(\cos\frac{t}{2}\right)^{\frac{1}{p}-1}
\int^{+\infty}_{0}\frac{y^{\alpha+\frac{1}{p}}}{(1+y^{2})^{\frac{\alpha+2}{2}}}
\frac{dy}{(1+y\cot\frac{t}{2})^{\alpha+\frac{2}{p}}}\\
&=2^{\alpha+\frac{2}{p}-1}c_{\alpha}\left(\sin\frac{t}{2}\right)^{\frac{1}{p}-1}
\left(\cos\frac{t}{2}\right)^{\frac{1}{p}-1}
\int^{+\infty}_{0}\frac{y^{\alpha+\frac{1}{p}}}{(1+y^{2})^{\frac{\alpha+2}{2}}}
\frac{dy}{(\sin\frac{t}{2}+y\cos\frac{t}{2})^{\alpha+\frac{2}{p}}}\\
&=2^{\alpha+\frac{1}{p}}c_{\alpha}(\sin t)^{\frac{1}{p}-1}
\int^{+\infty}_{0}\frac{y^{\alpha+\frac{1}{p}}}{(1+y^{2})^{\frac{\alpha+2}{2}}}
\frac{dy}{(\sin\frac{t}{2}+y\cos\frac{t}{2})^{\alpha+\frac{2}{p}}}.
\end{align*}

Let \begin{align*}G(t)=\int^{+\infty}_{0}\frac{y^{\alpha+\frac{1}{p}}}{(1+y^{2})^{\frac{\alpha+2}{2}}}
\frac{dy}{(\sin\frac{t}{2}+y\cos\frac{t}{2})^{\alpha+\frac{2}{p}}}, \quad t\in[0,\pi].
\end{align*}Changing variable with $y=\tan x,\, x\in[0, \frac{\pi}{2})$, we get
\begin{align*}
G(t)&=\int^{\frac{\pi}{2}}_{0}\frac{(\tan x)^{\alpha+\frac{1}{p}}}{(1+\tan^{2}x)^{\frac{\alpha+2}{2}}}
\frac{\sec^{2}xdx}{(\sin\frac{t}{2}+\tan x \cos\frac{t}{2})^{\alpha+\frac{2}{p}}}\\
&=\int^{\frac{\pi}{2}}_{0}\frac{(\tan x)^{\alpha+\frac{1}{p}}}{\sec^{\alpha}x}
\frac{(\cos x)^{\alpha+\frac{2}{p}}dx}{(\cos x\sin\frac{t}{2}+\sin x\cos\frac{t}{2})^{\alpha+\frac{2}{p}}}\\
&=\int^{\frac{\pi}{2}}_{0}\frac{(\sin x)^{\alpha+\frac{1}{p}}(\cos x)^{\alpha+\frac{1}{p}}}{\left(\sin
 (x+\frac{t}{2})\right)^{\alpha+\frac{2}{p}}}dx.
\end{align*}
Direct computation lead to
\begin{align*}
G''(t)=\frac{1}{4}(\alpha+\frac{2}{p})\int^{\frac{\pi}{2}}_{0}
\frac{(\sin x)^{\alpha+\frac{1}{p}}(\cos x)^{\alpha+\frac{1}{p}}}{\left(\sin
(x+\frac{t}{2})\right)^{2+\alpha+\frac{2}{p}}}\left((1+\alpha+\frac{2}{p})\cos^{2}(x+\frac{t}{2})+\sin^{2}(x+\frac{t}{2})\right)dx.
\end{align*}So, if $\alpha+\frac{2}{p}>0$, then it holds that $G''(t)>0$ for all $t\in [0,\pi]$. Furthermore, we have
\begin{align*}
G(0)=\int^{\frac{\pi}{2}}_{0}(\sin x)^{-\frac{1}{p}}(\cos x)^{\alpha+\frac{1}{p}}dx=\int^{\frac{\pi}{2}}_{0}
(\cos x)^{-\frac{1}{p}}(\sin x)^{\alpha+\frac{1}{p}}dx=G(\pi).
\end{align*}Recall that the Beta function has the following expression
 \begin{align*}
B(P,Q)=2\int^{\frac{\pi}{2}}_{0}\cos^{2P-1}t\sin^{2Q-1}tdt
\end{align*}for $P,\,Q>0$.
So,  \begin{align*}
G(0)=G(\pi)=\frac{1}{2}B(\frac{1+\alpha+\frac{1}{p}}{2},\frac{1-\frac{1}{p}}{2}).
\end{align*}Therefor, the convexity of $G(t)$ implies that
 \begin{align*}
G(t)\leq G(0)=G(\pi)=\frac{1}{2}B(\frac{1+\alpha+\frac{1}{p}}{2},\frac{1-\frac{1}{p}}{2})
\end{align*}for $t\in[0,\pi]$.

It follows from (\ref{3.14}) that
\begin{align*}
T^{*}((Th)^{p-1})&\leq\int^{1}_{-1}c_{\alpha}\frac{(1-r^{2})^{\alpha+\frac{1}{p}}}
{(1-2r\cos t+r^{2})^{\frac{\alpha+2}{2}}}dr
=2^{\alpha+\frac{1}{p}}c_{\alpha}(\sin t)^{\frac{1}{p}-1}G(t)\\
&\leq2^{\alpha+\frac{1}{p}-1}
c_{\alpha}(\sin t)^{\frac{1}{p}-1}B(\frac{1+\alpha+\frac{1}{p}}{2},\frac{1-\frac{1}{p}}{2}).
\end{align*}
Observe that
\begin{align*}
h^{p-1}(z)=2^{-\frac{p-1}{p}}(\sin t)^{-\frac{p-1}{p}}\left(\cos(\frac{\pi}{2p}-\frac{t}{p})\right)^{p-1}
\end{align*}for $t\in[0,\pi]$ with $z=e^{it}$.
Using the inequality $\cos^{p-1}(\frac{\pi}{2p}-\frac{t}{p})\geq\cos^{p-1}(\frac{\pi}{2p})$ for $t\in[0,\pi]$, we have
\begin{align*}
T^{*}((Th)^{p-1})\leq \tilde{C}(\alpha, p)h^{p-1}
\end{align*}with $\tilde{C}(\alpha, p)=\frac{2^{\alpha}c_{\alpha}B(\frac{1}{2}(1+\alpha+\frac{1}{p}),\frac{1}{2}(1-\frac{1}{p}))}
{\cos^{p-1}(\frac{\pi}{2p})}$. Therefore, by Lemma \ref{lem2.1}, we have
\begin{align*}
\int^{1}_{-1}|u(re^{is})|^{p}dr\leq   \tilde{C}(\alpha, p)\int^{2\pi}_{0}|u(e^{i\theta})|^{p}\frac{d\theta}{2\pi}:= C(\alpha, p)\int^{2\pi}_{0}|u(e^{i\theta})|^{p}d\theta.
\end{align*} with $ C(\alpha, p)=\frac{2^{\alpha-1}c_{\alpha}B(\frac{1}{2}(1+\alpha+\frac{1}{p}),\frac{1}{2}(1-\frac{1}{p}))}
{\pi\cos^{p-1}(\frac{\pi}{2p})}$.

If $\alpha=0$, then  recall the complement formula $\Gamma(x)\Gamma(1-x)=\frac{\pi}{\sin(x\pi)}$, $x\in(0,1)$, we have

\begin{align*}C(0, p)=\frac{c_{0}B(\frac{1}{2}(1+\frac{1}{p}),\frac{1}{2}(1-\frac{1}{p}))}
{2\pi\cos^{p-1}(\frac{\pi}{2p})}=\frac{\Gamma(\frac{1}{2}(1+\frac{1}{p})\Gamma(\frac{1}{2}(1-\frac{1}{p}))}
{2\pi\cos^{p-1}(\frac{\pi}{2p})}=\frac{1}
{2\sin(\frac{1}{2}(1-\frac{1}{p})\pi)\cos^{p-1}(\frac{\pi}{2p})}=\frac{1}
{2}\sec^{p}(\frac{\pi}{2p}).\end{align*}
Then  it reduces to the case of harmonic functions.
\end{proof}

\medskip

{\bf Declarations}\\

{\bf Conflict of interests}\quad The authors declare that they have no conflict of interest.\\

{\bf Data availability statement}\quad This manuscript has no associated date.

\medskip

\end{document}